%% file: qr-practice.tex
\documentclass[final]{siamltex}

\usepackage{color}

\usepackage{amsmath,amssymb}
\usepackage{algorithm,algpseudocode}  
\usepackage{listings,url,verbatim}    
\usepackage{multicol}
\usepackage{multirow}         
\usepackage[pdftex]{graphicx} 
\usepackage{subfigure}        

\graphicspath{{./FIGURES/}}

\title{Implementing Communication-Optimal Parallel and Sequential QR
  Factorizations}

\author{James Demmel\thanks{Computer Science Division and Mathematics
    Department, UC Berkeley, CA 94720-1776, USA ({\tt
      demmel@cs.berkeley.edu}).}  \and Laura Grigori\thanks{INRIA
    Saclay-Ile de France, Laboratoire de Recherche en Informatique,
    Bat 490 University Paris-Sud 11, 91405 Orsay, France ({\tt
      laura.grigori@inria.fr}).}  \and Mark Hoemmen\thanks{Computer
    Science Division, UC Berkeley, CA 94720-1776, USA ({\tt
      mhoemmen@eecs.berkeley.edu}).}  \and Julien
  Langou\thanks{Department of Mathematical and Statistical Sciences, University of
    Colorado Denver, CO 80202, USA ({\tt julien.langou@ucdenver.edu}).}}

\begin{document}
\maketitle

\lstset{basicstyle=\small,showstringspaces=false,language=C} 

\begin{abstract}
  We present parallel and sequential dense QR factorization algorithms
  for tall and skinny matrices and general rectangular matrices that both
  minimize communication, and are as \emph{stable} as Householder QR.
  The sequential and parallel algorithms for tall and skinny matrices lead
  to significant speedups in practice over some of the existing
  algorithms, including LAPACK and ScaLAPACK, for example up to 6.7x 
  over ScaLAPACK.
  The parallel algorithm
  for general rectangular matrices is estimated to show significant
  speedups over ScaLAPACK, up to 22x over ScaLAPACK.
\end{abstract}



\pagestyle{myheadings}
\thispagestyle{plain}
\markboth{J. DEMMEL, L. GRIGORI, M. HOEMMEN, AND J. LANGOU}{Implementing communication-optimal QR factorization}

\input{pr-introduction}

\input{pr-tsqr-algebra}

\input{pr-tsqr-opt}

\input{pr-tsqr-impl}

\input{pr-tsqr-other}

\input{pr-caqr}

\input{pr-experiments}

\input{pr-conclusions}

\bibliographystyle{siam}
\bibliography{../TechReport2007/qr}
\end{document}

%% file: pr-introduction.tex
\section{Introduction}

In this paper we present parallel and sequential dense QR
factorization algorithms that both minimize communication, and are as
\emph{stable} as Householder QR. (That is to say normwise backward
stable.)  Communication refers to messages that are sent over a
network in the parallel case, and to data movement between different
levels of memory hierarchy in the sequential case.  The first set of
algorithms, ``Tall Skinny QR'' (TSQR), are for matrices for which the
number of rows is much larger than the number of columns, and which
have their rows distributed over processors in a one-dimensional (1-D)
block row layout.  The second set of algorithms,
``Communication-Avoiding QR'' (CAQR), are for general rectangular
matrices distributed using a two-dimensional (2-D) block cyclic
layout.  Some of these algorithms are new, and some are based on
existing work.

The new algorithms are superior in both theory and practice.
In~\cite{THEORY} we show that the new sequential and parallel 
algorithms, for both 1-D layout TSQR and 2-D block cyclic layout CAQR,
are communication optimal (modulo polylogarithmic factors), that is 
they are optimal in bandwidth and latency costs. This assumes $O(n^3)$
algorithms (non-Strassen-like).
We also observe in~\cite{THEORY} that
LAPACK's corresponding sequential factorization and ScaLAPACK's
parallel QR factorization perform asymptotically more communication.
The sequential recursive QR factorization algorithm by Elmroth and
Gustavson \cite{elmroth2000applying} attains the lower bound on the
volume of data transferred in some special cases, though not the lower
bound on the number of block transfers.

In this paper, we focus on the implementation and performance results
of these algorithms.  The efficient implementation of the QR
factorizations of tall and skinny matrices distributed in a 1-D layout
is very important, since this operation arises in a wide range of
applications.  We cite three important examples.  Block iterative
methods frequently compute the QR factorization of a tall and skinny
dense matrix.  This includes algorithms for solving linear systems $Ax
= B$ with multiple right-hand sides (such as variants of GMRES, QMR,
or CG \cite{vital:phdthesis:90,Freund:1997:BQA,oleary:80}), as well as
block iterative eigensolvers (for a summary of such methods, see
\cite{templatesEigenBai,templatesEigenLehoucq}).  Many of these
methods have widely used implementations, on which a large community
of scientists and engineers depends for their computational tasks.
Examples include TRLAN (Thick Restart Lanczos), BLZPACK (Block
Lanczos), Anasazi (various block methods), and PRIMME (block
Jacobi-Davidson methods)
\cite{TRLANwebpage,BLZPACKwebpage,BLOPEXwebpage,irbleigs,TRILINOSwebpage,PRIMMEwebpage}.
Eigenvalue computation is particularly sensitive to the accuracy of
the orthogonalization; two recent papers suggest that large-scale
eigenvalue applications require a stable QR factorization
\cite{lehoucqORTH,andrewORTH}. Our approach is based on Householder
reflections so it is unconditionnally normwise backward stable, as
opposed to other standard used methods as Gram-Schmidt or CholeskyQR
(see Section~\ref{pr:TSQR:other}).

Recent research has reawakened an interest in alternate formulations
of Krylov subspace methods, called \emph{$s$-step Krylov methods}, in
which some number $s$ steps of the algorithm are performed all at
once, in order to reduce communication.  Demmel et al.\ review the
existing literature and discuss new advances in this area
\cite{demmel2008comm}.  Such a method generates some basis for the
Krylov subspace, and then a QR factorization is used to orthogonalize
the basis vectors.  This is an ideal application for TSQR, and in fact
inspired its (re-)discovery.
 
Householder QR decompositions of tall and skinny matrices also
comprise the panel factorization step for typical QR factorizations of
matrices in a more general, two-dimensional layout.  This includes the
current parallel QR factorization routine \lstinline!PDGEQRF! in
ScaLAPACK, as well as ScaLAPACK's out-of-DRAM QR factorization
\lstinline!PFDGEQRF!~\cite{dazevedo1997design_paper}.
Both algorithms use a standard column-based
Householder QR for the panel factorizations, but in the parallel case
this is a latency bottleneck, and in the out-of-DRAM case it is a
bandwidth bottleneck.  Our CAQR algorithm for computing the QR
factorization of general rectangular matrices uses TSQR for its panel
factorization.  That's how CAQR removes the latency bottleneck in the
parallel case and the bandwidth bottleneck in the sequential case.

The main insight behind TSQR algorithm is to perform the QR
factorization of a tall-skinny matrix as a reduction operation.  This
idea itself is not novel (see for example,
\cite{bggl:06,buttari2007class,cunha2002new,golub1988parallel,gunter2005parallel,kurzak2008qr,pothen1989distributed,quintana-orti2008scheduling,rabani2001outcore}),
but we have a number of optimizations and generalizations:
\begin{itemize}
\item Our TSQR algorithm can compute most of its floating-point
  operations by using the best available sequential QR factorization.
  In particular, we can achieve significant speedups by invoking
  Elmroth and Gustavson's recursive QR (see
  \cite{elmroth1998new,elmroth2000applying}).
\item We use TSQR as a building block for CAQR, the parallel
  factorization of rectangular matrices with a two-dimensional block
  cyclic layout.  To our knowledge, parallel CAQR is a new algorithm.
\item We explain how TSQR can work on general reduction trees.  This
  flexibility lets schedulers overlap communication and computation,
  and minimize communication for more complicated and realistic
  computers with multiple levels of parallelism and memory hierarchy
  (e.g., a system with disk, DRAM, and cache on multiple boards each
  containing one or more multicore chips of different clock speeds,
  along with compute accelerator hardware like GPUs).
\end{itemize}
In practice, parallel TSQR leads to significant speedups on several
machines:
\begin{itemize}
\item up to $6.7\times$ on 16 processors of a Pentium III cluster, for
  a $100,000 \times 200$ matrix; and
\item up to $4\times$ on 32 processors of a BlueGene/L, for a
  $1,000,000 \times 50$ matrix.
\end{itemize}
Some of this speedup is enabled by TSQR being able to use a much
better local QR decomposition than ScaLAPACK can use, such as the
recursive variant by Elmroth and Gustavson \cite{elmroth2000applying}.
We have also implemented sequential TSQR on a laptop for matrices that
do not fit in DRAM, so that slow memory is disk.  This requires a
special implementation in order to run at all, since virtual memory
does not accommodate matrices of the sizes we tried.  By extrapolating
runtime from matrices that do fit in DRAM, we can say that our
out-of-DRAM implementation was as little as $2\times$ slower than the
predicted runtime as though DRAM were infinite.

We also estimate the performance of our parallel CAQR algorithm (whose
actual implementation and measurement is current work), yielding
predicted speedups over ScaLAPACK's \lstinline!PDGEQRF! of up to
$22.9\times$ on a model Petascale machine.  The best speedups occur
for the largest number of processors used, and for matrices that do
not fill all of memory, since in this case latency costs dominate. In
general, when the largest possible matrices are used, computation
costs dominate the communication costs and improved communication does
not help.

The rest of the paper is organized as follows.
Section~\ref{pr:S:algebra} describes the algebra of the TSQR algorithm
and shows how the parallel and sequential versions correspond to
different trees.  Section~\ref{pr:S:localQR} shows that the local QR
decompositions in TSQR can be further optimized by exploiting the structure of the matrices involved.  We also explain how
to apply the $Q$ factor from TSQR efficiently, which is needed both
for CAQR and other applications.  Section \ref{S:TSQR:impl}
describes the parallel and sequential TSQR algorithms, and presents a
performance model for each.  Next, Section~\ref{pr:TSQR:other}
describes other ``tall skinny QR'' algorithms, such as CholeskyQR and
Gram-Schmidt, and compares their cost and numerical stability to that
of TSQR.  This section shows that TSQR is the only algorithm that
simultaneously minimizes communication and is numerically stable.
Section~\ref{pr:S:CAQR} describes the parallel CAQR algorithm and
constructs a performance model.  Section~\ref{pr:experiments} presents
the platforms used for testing, and discusses the TSQR and CAQR
performance results. Section~\ref{pr:S:concl} concludes the paper.

%% file: pr-tsqr-algebra.tex
\section{TSQR matrix algebra}
\label{pr:S:algebra}

In this section, we illustrate the insight behind the TSQR algorithm
for computing the QR factorization of an $m \times n$ matrix $A$
partitioned in $P$ block rows, that is $A = [A_0 ; A_1 ; \cdots ;
A_{P-1}]$.  We use Matlab notation, so that the $A_i$ are stacked atop
one another.

TSQR defines a family of algorithms, in which the QR factorization of
$A$ is obtained by performing a sequence of QR factorizations until
the lower trapezoidal part of $A$ is annihilated and the final $R$
factor is obtained.  The QR factorizations are performed on block rows
of $A$ and on previously obtained $R$ factors, stacked atop one
another.  We call the pattern followed during this sequence of QR
factorizations a reduction tree.  We begin with parallel TSQR, for
which the reduction tree is a binary tree, and later show sequential
TSQR on a linear tree.  We consider the simple example of $P=4$.

Parallel TSQR starts with the independent computation of the QR
factorization of each block row:
\[
A =
\begin{pmatrix}
A_0 \\
A_1 \\
A_2 \\
A_3 \\
\end{pmatrix}
=
\begin{pmatrix}
  Q_{00} R_{00} \\
  Q_{10} R_{10} \\
  Q_{20} R_{20} \\
  Q_{30} R_{30} \\
\end{pmatrix}.
\]
This is ``stage 0'' of the computation, hence the second subscript 0
of the $Q$ and $R$ factors.  The first subscript indicates the block
index at that stage.  Stage 0 operates on the $P = 4$ leaves of the
tree.  After this stage, there are $P = 4$ of the $R$ factors.  We
group them into successive pairs $R_{i,0}$ and $R_{i+1,0}$, and do the
QR factorizations of grouped pairs in parallel:
\[
\begin{pmatrix}
  R_{00} \\
  R_{10} \\ \hline
  R_{20} \\ 
  R_{30} \\
\end{pmatrix}
= 
\begin{pmatrix}
  \begin{pmatrix} 
    R_{00} \\
    R_{10} \\
  \end{pmatrix} \\ \hline
  \begin{pmatrix}
    R_{20} \\
    R_{30} \\
  \end{pmatrix}
\end{pmatrix}
=
\begin{pmatrix}
  Q_{01} R_{01} \\ \hline
  Q_{11} R_{11} \\ 
\end{pmatrix}.
\]
This is stage 1, as the second subscript of the $Q$ and $R$ factors
indicates.  We iteratively perform stages until there is only one $R$
factor left, which is the root of the tree:
\[
\begin{pmatrix}
  R_{01} \\
  R_{11} \\
\end{pmatrix}
=
Q_{02} R_{02}.
\]

If we were to compute all the above $Q$ factors explicitly as square
matrices, which we do not, 
each of the $Q_{i0}$ would be $m/P \times m/P$, and $Q_{ij}$
for $j > 0$ would be $2n \times 2n$.  The final $R_{02}$ factor would
be $m \times n$ upper triangular (or $n \times n$ upper triangular in
a ``thin QR'' factorization).

Equation \eqref{eq:TSQR:algebra:par4:final} shows the whole
factorization:
\begin{equation}\label{eq:TSQR:algebra:par4:final}
A = 
\left(
  \begin{array}{c | c | c | c}
    Q_{00} & & & \\ \hline
    & Q_{10} & & \\ \hline
    & & Q_{20} & \\ \hline
    & & & Q_{30} \\ 
  \end{array}
\right)
\cdot
\left(
  \begin{array}{c | c}
    \tilde{Q}_{01} &       \\ \hline
          & \tilde{Q}_{11} \\ 
  \end{array}
\right)
\cdot
\tilde{Q}_{02} \cdot R_{02},
\end{equation}
in which $\tilde{Q}_{ij}$ with $j > 1$ are the matrices $Q_{ij}$
extended by the identity to match the dimensions $m \times m$ of the first
$Q$ factor.

The product of the first three matrices is orthonormal, because each
of these three matrices is.  Since the QR decomposition is essentially
unique (it is unique modulo signs of diagonal entries of $R_{02}$),
this is the QR decomposition of $A$ and $R_{02}$ is the $R$ factor of
$A$.

Note the binary tree structure in the nested pairs of $R$ factors.
This tree structure and the underlined TSQR algorithm can be
visualized using a similar notation to~\cite{THEORY}:

\begin{center}
\setlength{\unitlength}{.5cm}
\begin{picture}(7,4)

\put(0.5,0.5){$A_3$}
\put(0.5,1.5){$A_2$}
\put(0.5,2.5){$A_1$}
\put(0.5,3.5){$A_0$}

\put(1.5,0.5){$\rightarrow$}
\put(1.5,1.5){$\rightarrow$}
\put(1.5,2.5){$\rightarrow$}
\put(1.5,3.5){$\rightarrow$}

\put(2.5,0.5){$R_{30}$}
\put(2.5,1.5){$R_{20}$}
\put(2.5,2.5){$R_{10}$}
\put(2.5,3.5){$R_{00}$}

\put(3.5,0.65){$\nearrow$}
\put(3.5,1.35){$\searrow$}
\put(3.5,2.65){$\nearrow$}
\put(3.5,3.35){$\searrow$}

\put(4.5,1.0){$R_{11}$}
\put(4.5,3.0){$R_{01}$}

\put(5.6,1.5){$\nearrow$}
\put(5.6,2.5){$\searrow$}

\put(6.5,2.0){$R_{02}$}

\end{picture}
\end{center}
where the arrows pointing to an $R$ factor highlight the matrices
whose QR factorization, when stacked atop one another, results in this
$R$ factor.  This representation illustrates well the parallelism
available in the algorithm as well.  The $R$ nodes of the tree
represent local QR factorizations, that is computations performed by
one processor, and the arrows between $R$ factors represent
communication.

Sequential TSQR uses a similar factorization process, but with a
``flat tree'' (a linear chain).  We start with the same block row
decomposition as with parallel TSQR, but begin with a QR factorization
of $A_0$, rather than of all the block rows:
\[
A =
\begin{pmatrix}
A_0 \\
A_1 \\
A_2 \\
A_3 \\
\end{pmatrix}
=
\begin{pmatrix}
Q_{00} R_{00} \\
A_1 \\
A_2 \\
A_3 \\
\end{pmatrix}.
\]
This is ``stage 0'' of the computation, hence the second subscript 0
of the $Q$ and $R$ factor.  We retain the first subscript for
generality, though in this example it is always zero.  We then combine
$R_{00}$ and $A_1$ using a QR factorization:
\[
\begin{pmatrix}
R_{00} \\
A_1 \\
A_2 \\
A_3 \\
\end{pmatrix}
=
\begin{pmatrix}
R_{00} \\
A_1 \\ \hline
A_2 \\
A_3 \\
\end{pmatrix}
=
\begin{pmatrix}
  Q_{01} R_{01} \\ \hline
  A_2 \\
  A_3 \\
\end{pmatrix}
\]
We continue this process until we run out of $A_i$ factors.  Here, the
$A_i$ blocks are $m/P \times n$.  If we were to compute all the above
$Q$ factors explicitly as square matrices, which we do not,
then $Q_{00}$ would be $m/P
\times m/P$ and $Q_{0j}$ for $j > 0$ would be $2m/P \times 2m/P$.  The
final $R$ factor, as in the parallel case, would be $m \times n$ upper
triangular (or $n \times n$ upper triangular in a ``thin QR'').

The resulting factorization has the following structure:
\begin{equation}\label{eq:TSQR:algebra:seq4:final}
A 
=
\left(
  \begin{array}{c | c | c | c}
    Q_{00} & & & \\ \hline
    & I & & \\ \hline
    & & I & \\ \hline
    & & & I \\ 
  \end{array}
\right)
\cdot
\left(
\begin{array}{c | c | c}
  \tilde{Q}_{01} &   &    \\ \hline
        & I &    \\ \hline
        &   & I  \\
\end{array}
\right)
\cdot
\left(
  \begin{array}{ c | c}
     \tilde{Q}_{02} & \\ \hline
     & I \\
  \end{array}
\right)
\cdot
\tilde{Q}_{03}
\cdot
R_{03}.
\end{equation}
where $\tilde{Q}_{0j}$ with $j > 1$ are the matrices $Q_{0j}$ extended
by the identity to match the dimensions of the equation.  The above $I$
factors are $m/P \times m/P$.

Sequential TSQR and the flat tree structure on which the factorization
executes is illustrated using the ``arrow'' abbreviation as:

\begin{center}
\setlength{\unitlength}{.5cm}
\begin{picture}(7,4)

\put(0.5,0.5){$A_3$}
\put(0.5,1.5){$A_2$}
\put(0.5,2.5){$A_1$}
\put(0.5,3.5){$A_0$}

\put(1.5,1.0){\vector(3,1){7}}
\put(1.5,1.75){\vector(3,1){5}}
\put(1.5,2.75){\vector(4,1){3}}
\put(1.5,3.75){\vector(1,0){1}}

\put(2.5,3.5){$R_{00}$}
\put(3.5,3.75){\vector(1,0){1}}
\put(4.5,3.5){$R_{01}$}
\put(5.55,3.75){\vector(1,0){.8}}
\put(6.5,3.5){$R_{02}$}
\put(8.0,3.75){\vector(1,0){.5}}
\put(8.5,3.5){$R_{03}$}

\end{picture}
\end{center}

A similar algorithm, but with a bottom-up traversal of the flat tree,
can also be formulated.  The flat-tree approach is similar to the
updating techniques proposed for out-of-core
computations~\cite{bggl:06,gunter2005parallel}, or for
multicore~\cite{buttari2007class,quintana-orti2008scheduling} and Cell
processors~\cite{kurzak2008qr}.






The sequential algorithm differs from the parallel one in that it does
not factor the individual blocks of the input matrix $A$, excepting
$A_0$.  This is because in the sequential case, a bit more than one
block of $A$ can be loaded into working memory.  In the fully parallel
case, each block of $A$ resides in some processor's working memory.
It then pays to factor all the blocks before combining them, as this
reduces the volume of communication (only the triangular $R$ factors
need to be exchanged) and reduces the amount of arithmetic performed
at the next level of the tree.  In contrast, the sequential algorithm
never writes out the intermediate $R$ factors, so it does not need to
convert the individual $A_i$ into upper triangular factors.

The above two algorithms are extreme points in a large set of possible
QR factorization methods, parametrized by the tree structure.  Our
version of TSQR is novel because it works on any tree.  In general,
the optimal tree may depend on both the architecture and the matrix
dimensions.  This is because TSQR is a reduction (as we will discuss
further in Section \ref{pr:SS:tsqr_red}).  Trees of types other than
binary often result in better reduction performance, depending on the
architecture (see e.g., \cite{nishtala2008performance}).  Throughout
this paper, we discuss two examples -- the binary tree and the flat
tree -- as easy extremes for illustration.  It is shown
in~\cite{THEORY} that the binary tree
is optimal in
the number of stages
and messages in the parallel case, and that the flat tree
is optimal in
the number and volume of input matrix reads and writes in the
sequential case.  Methods for finding the best tree in general
are future work.  We expect to use a non-binary tree in the case
of real-world systems with multiple levels of memory hierarchy and
multiple, possibly heterogeneous processors, although in this paper we
do not address this issue.

%% file: pr-tsqr-opt.tex
\section{Optimizations for TSQR}
\label{pr:S:localQR}

Although TSQR achieves its performance gains because it optimizes
communication, the local QR factorizations lie along the critical path
of the algorithm.  The parallel cluster benchmark results in Section
\ref{pr:experiments} show that optimizing the local QR factorizations
can improve performance significantly.  In this section, we outline a
few of these optimizations, and hint at how they affect the
formulation of the general CAQR algorithm in Section \ref{pr:S:CAQR}.

\subsection{Optimizing local factorizations in TSQR}
\label{pr:TSQR_localopt} 

Most of the QR factorizations performed during TSQR involve matrices
consisting of one or more triangular factors stacked atop one another.
We can ignore this zero structure and still get a correct
factorization, but if we do we will do several times as many floating
point operations as necessary (up to 5$\times$ in the parallel case
and 2$\times$ in the sequential case). Previous authors have suggested
using Givens rotations to avoid this \cite{pothen1989distributed}, but
this would make it hard to achieve Level 3 BLAS performance.

Our observation is that not only it is possible to use blocked
Householder transformations that both do minimal arithmetic and permit
Level 3 BLAS performance, but in fact we can organize the algorithm to get
{\em better} Level 3 BLAS performance than conventional QR decomposition.
The empirical data justifying this claim appears in
Section~\ref{pr:experiments}, but we outline the idea here.

We illustrate with the QR decomposition of a pair $[R_0 ; R_1 ]$ of
5-by-5 triangular matrices.  Their sparsity pattern, and that of the
Householder vectors from their QR decomposition are shown below:

\begin{equation}\label{eq:fact_2rs}
\begin{pmatrix}
R_0 \\
R_1\\
\end{pmatrix}
=
\begin{pmatrix}
x & x & x & x & x \\
  & x & x & x & x \\
  &   & x & x & x \\
  &   &   & x & x \\
  &   &   &   & x \\
x & x & x & x & x \\
  & x & x & x & x \\
  &   & x & x & x \\
  &   &   & x & x \\
  &   &   &   & x \\
\end{pmatrix}
\Longrightarrow
{\rm Householder} = 
\begin{pmatrix}
1 &   &   &   &   \\
  & 1 &   &   &   \\
  &   & 1 &   &   \\
  &   &   & 1 &   \\
  &   &   &   & 1 \\
x & x & x & x & x  \\
  & x & x & x & x \\
  &   & x & x & x \\
  &   &   & x & x \\
  &   &   &   & x \\
\end{pmatrix}.
\end{equation}

This picture suggests that it is straightforward to adapt both the
unblocked Householder decomposition and its blocked version in
\cite{schreiber1989storage}, by storing the Householder vectors on top
of the zeroed-out entries as usual, and simply by changing the lengths
of the vectors involved in updates of the trailing matrix.  For the
case of two $n$-by-$n$ triangular matrices, exploiting this structure
lowers the operation count to $\frac{2}{3}n^3$ from about
$\frac{10}{3}n^3$.  It is also possible to do this when $q$ triangles
are stacked atop one another, or when a triangle is stacked atop a
rectangular block as in sequential TSQR.  Most importantly, we can
apply Elmroth and Gustavson's recursive QR algorithm
\cite{elmroth2000applying} to the matrices in fast memory (in the
sequential case) or local processor memory (in the parallel case).

\subsection{Trailing matrix update}\label{SS:TSQR:localQR:trailing}

Section \ref{pr:S:CAQR} will describe how to use TSQR to factor
matrices in general 2-D layouts.  For these layouts, once the current
panel (block column) has been factored, the panels to the right of the
current panel cannot be factored until the transpose of the current
panel's $Q$ factor has been applied to them.  This is called a
\emph{trailing matrix update} and consists of a sequence of
applications of local $Q^T$ factors to groups of trailing matrix
blocks. The update lies along the critical path of the algorithm, and
consumes most of the floating-point operations in general.  We now
explain how to do one of these local $Q^T$ applications.

Let the number of rows in a block be $2n$, and the number of columns
in a block be $n$.  Suppose that we want to apply the local $Q^T$
factor from the above $2n \times n$ matrix factorization, to two
blocks $C_0$ and $C_1$ of a trailing matrix panel.  
$C_0$ and $C_1$ may have more than $n$ columns.
Our goal is to perform the operation
\[
\begin{pmatrix}
 R_0 & C_0 \\
 R_1 & C_1 \\
\end{pmatrix}
= 
\begin{pmatrix}
QR & C_0 \\
   & C_1 \\
\end{pmatrix}
=
Q \cdot
\begin{pmatrix}
R & \hat{C}_0 \\
  & \hat{C}_1 \\
\end{pmatrix},
\]
in which $Q$ is the local $Q$ factor and $R$ is the local $R$ factor
of $[R_0; R_1]$.  When the YT representation is used for
$Q$~\cite{schreiber1989storage}, the update of the trailing matrices
takes the following form:
\[
\begin{pmatrix}
\hat{C_0} \\
\hat{C_1} \\
\end{pmatrix}
:=
\begin{pmatrix}
I -
\begin{pmatrix}
 I \\
Y_1 \\
\end{pmatrix}
\cdot 
\begin{array}{cc} 	
T^T \\
\\
\end{array}
\cdot
\begin{pmatrix}
 I \\
Y_1 \\
\end{pmatrix}^T
\end{pmatrix}
\begin{pmatrix}
C_0 \\
C_1 \\
\end{pmatrix}.
\]
If we let
\[
D := C_0 + Y_1^T C_1 
\]
be the ``inner product'' part of the update operation formulas, then
we can rewrite the update formulas as
\[
\begin{aligned}
\hat{C_0} &:= C_0 - T^T D, \\
\hat{C_1} &:= C_1 - Y_1 T^T D, \\
\end{aligned}
\]

In a parallel algorithm, there are many different ways to perform this
update.  The data dependencies impose a \emph{directed acyclic graph} (DAG)
on the flow of data between processors.  One can find the best way to
do the update by realizing an optimal computation schedule on the DAG.
In Section~\ref{pr:S:CAQR} we will see a straightforward 
schedule of this computation.

%% file: pr-tsqr-impl.tex
\section{Parallel and sequential TSQR}
\label{S:TSQR:impl}

In this section, we describe the TSQR factorization algorithm in more
detail.  We also build a performance model of the algorithm, based on
a simple machine model.  We predict floating-point performance by
counting floating-point operations and multiplying them by $\gamma$,
the inverse peak floating-point performance.  We use the
``alpha-beta'' or latency-bandwidth model of communication, in which a
message of size $n$ floating-point words takes time $\alpha + \beta n$
seconds.  The $\alpha$ term represents message latency (seconds per
message), and the $\beta$ factor inverse bandwidth (seconds per
floating-point word communicated).  We also apply the alpha-beta model
to communication between levels of the memory hierarchy in the
sequential case.  We restrict our model to describe only two levels at
one time: fast memory (which is smaller) and slow memory (which is
larger).

Parallel TSQR performs $2mn^2/P + \frac{2n^3}{3}\log P$ flops,
compared to the $2mn^2/P - 2n^3/(3P)$ flops performed by ScaLAPACK's
parallel QR factorization \lstinline!PDGEQRF!, but requires $2n$ times
fewer messages.  The sequential TSQR factorization performs the same
number of flops as sequential blocked Householder QR, but requires
$\mathcal{O}(n)$ times fewer transfers between slow and fast memory, and a
factor of $\mathcal{O}(mn/W)$ times fewer words transferred, in which
$W$ is the fast memory size. Note that $mn/W$ is how many times larger
the matrix is than the fast memory.

\subsection{TSQR as a reduction}
\label{pr:SS:tsqr_red}

Section~\ref{pr:S:algebra} explained the algebra of the TSQR
factorization.  It outlined how to reorganize the parallel QR
factorization as a tree-structured computation, in which groups of
neighboring processors combine their $R$ factors, perform (possibly
redundant) QR factorizations, and continue the process by
communicating their $R$ factors to the next set of neighbors.
Sequential TSQR works in a similar way, except that communication
consists of moving matrix factors between slow and fast memory.  This
tree structure uses the same pattern of communication found in a
reduction or all-reduction.  We can say TSQR factorization is itself
an (all-)reduction, in which additional data (the components of the
$Q$ factor) is stored at each node of the (all-)reduction tree.
Applying the $Q^T$ factor is also a(n) (all-)reduction; while applying the $Q$
factor is a broadcast-like algorithm.

However, TSQR has requirements that differ from the standard (all-)reduction.
For example, if the $Q$ factor is desired, then TSQR must store intermediate
results (the local $Q$ factor from each level's computation with neighbors) at
interior nodes of the tree.  This requires reifying and preserving the
(all-)reduction tree for later invocation by users.  Typical (all-)reduction
interfaces, such as those provided by MPI or OpenMP, do not easily allow this
(see e.g., \cite{gropp1999using}).


When TSQR is implemented with an all-reduction, then the final $R$
factor is replicated over all the processors.  This is especially
useful for Krylov subspace methods.  If TSQR is implemented with a
simple reduction, then the final $R$ factor is stored only on one
processor.  This avoids redundant computation, and is useful both for
block column factorizations for 2-D block (cyclic) matrix layouts, and
for solving least squares problems when the $Q$ factor is not needed.

\subsection{Factorization}

We now describe the parallel and sequential TSQR factorizations for
the 1-D block row layout.  (We omit the obvious generalization to a
1-D block cyclic row layout.)

Parallel TSQR computes an $R$ factor which is duplicated over all the
processors, and a $Q$ factor which is stored implicitly in a
distributed way.  The algorithm overwrites the lower trapezoid of
$A_{i}$ with the set of Householder reflectors for that block, and the
$\tau$ array of scaling factors for these reflectors is stored
separately.  The matrix $R_{i,k}$ is stored as an $n \times n$ upper
triangular matrix for all stages $k$.  Algorithm \ref{Alg:TSQR:par}
shows an implementation of parallel TSQR, based on an all-reduction.

\begin{algorithm}[h]
\caption{Parallel TSQR}
\label{Alg:TSQR:allred:blkrow}
\label{Alg:TSQR:par}
\begin{algorithmic}[1]
\Require{$\Pi$ is the set of $P$ processors}
\Require{All-reduction tree with height $L$.  If $P$ is a power of two
  and we want a binary all-reduction tree, then $L = \log_2 P$.}
\Require{$i \in \Pi$: my processor's index}
\Require{The $m \times n$ input matrix $A$ is distributed in a 1-D 
  block row layout over the processors; $A_{i}$ is the block of rows 
  belonging to processor $i$}.
\State{Compute $[Q_{i,0}, R_{i,0}] := qr(A_{i})$ using sequential 
  Householder QR}\label{Alg:TSQR:allred:blkrow:QR1}
\For{$k$ from 1 to $L$}
  \If{I have any neighbors in the all-reduction tree at this level}
    \State{Send (non-blocking) $R_{i,k-1}$ to each neighbor not myself}
    \State{Receive (non-blocking) $R_{j,k-1}$ from each neighbor $j$ not myself}
    \State{Wait until the above sends and receives
      complete}\Comment{Note: \emph{not} a global barrier.}
    \State{Stack the upper triangular $R_{j,k-1}$ from all neighbors 
      (including my own $R_{i,k-1}$), by order of processor ids, into 
      a $qn \times n$ array $C$, in which $q$ is the number of 
      neighbors.}
    \State{Compute $[Q_{i,k}, R_{i,k}] := qr(C)$}
  \Else
    \State{$R_{i,k} := R_{i,k-1}$}
    \State{$Q_{i,k} := I_{n \times n}$}\Comment{Stored implicitly}
  \EndIf
  \State{Processor $i$ has an implicit representation of its block column
    of $Q_{i,k}$.  The blocks in the block column are $n \times n$
    each and there are as many of them as there are neighbors at stage
    $k$ (including $i$ itself).  We don't need to compute the blocks
    explicitly here.}
\EndFor
\Ensure{$R_{i,L}$ is the $R$ factor of $A$, for all processors $i \in \Pi$.}
\Ensure{The $Q$ factor is implicitly represented by $\{Q_{i,k}\}$:
  $i \in \Pi$, $k \in \{0, 1, \dots, L\}\}$.}
\end{algorithmic}
\end{algorithm}

At the leaf nodes of the TSQR tree (step 1 of TSQR algorithm), each
processor computes a QR factorization of an $m/P \times n$ matrix.
This factorization involves around $2n^2 m / P - 2 n^3 / 3$ flops.
For all the other nodes of the TSQR tree (step 2 of the TSQR
Algorithm), two processors perform redundantly the QR factorization of
a $2n \times n$ matrix formed by two upper triangular matrices.  The
number of flops performed on the critical path of TSQR is $2n^2 m / P
- 2 n^3 / 3 + 2 n^3 / 3 \log P$.  Thus, the run time of the TSQR
algorithm is estimated to be

\begin{equation}
\label{eqn:TSQR_par_runtime}
\text{Time}_{\text{Par.\ TSQR}}(m,n,P) = 
\left(
  \frac{2mn^2}{P} + \frac{2n^3}{3} \log P
\right) \gamma + 
\left(
  \frac{1}{2} n^2 \log P
\right) \beta + 
\left( \log P \right) \alpha 
\; \; .
\end{equation}

Sequential TSQR begins with an $m \times n$ matrix $A$ stored in slow
memory.  The matrix $A$ is divided into $P$ blocks $A_0$, $A_1$,
$\dots$, $A_{P-1}$, each of size $m/P \times n$.  (Here, $P$ has
nothing to do with the number of processors; it is chosen to
minimize latency, i.e. as small as possible subject to the memory
constraint described below.)  Each block of $A$ is
loaded into fast memory in turn, combined with the $R$ factor from the
previous step using a QR factorization, and the resulting $Q$ factor
written back to slow memory.  Thus, only one $m/P \times n$ block of
$A$ resides in fast memory at one time, along with an $n \times n$
upper triangular $R$ factor.  
Thus we choose $P$ as small as possible subject to the memory
constraint $\frac{mn}{P} + \frac{n(n+1)}{2} \leq W$.
Sequential TSQR computes an $n \times n$
$R$ factor which ends up in fast memory, and a $Q$ factor which is
stored implicitly in slow memory as a set of blocks of Householder
reflectors.  Algorithm \ref{Alg:TSQR:seq} shows an implementation of
sequential TSQR.

\begin{algorithm}[h]
\caption{Sequential TSQR}\label{Alg:TSQR:seq}
\begin{algorithmic}[1]
\Require{The $m \times n$ input matrix $A$, stored in slow memory, 
  is divided into $P$ row blocks $A_0$, $A_1$, $\dots$, $A_{P-1}$}
\State{Load $A_0$ into fast memory}
\State{Compute $[Q_{00}, R_{00}] := qr(A_{0})$ using standard
  sequential QR.   Here, the $Q$ factor is represented implicitly by
  an $m/P \times n$ lower triangular array of Householder reflectors 
  $Y_{00}$ and their $n$ associated scaling factors $\tau_{00}$}
\State{Write $Y_{00}$ and $\tau_{00}$ back to slow memory; keep
  $R_{00}$ in fast memory}
\For{$k = 1$ to $P - 1$}
    \State{Load $A_k$}
    \State{Compute $[Q_{01}, R_{01}] = qr([R_{0,k-1}; A_k])$.   Here, the $Q$ factor 
      is represented implicitly by a full $m/P \times n$ array of 
      Householder reflectors $Y_{0k}$ and their $n$ associated 
      scaling factors $\tau_{0k}$.}
    \State{Write $Y_{0k}$ and $\tau_{0k}$ back to slow memory; 
      keep $R_{0k}$ in fast memory}
\EndFor
\Ensure{$R_{0,P-1}$ is the $R$ factor in the QR factorization of $A$,
  and is in fast memory}
\Ensure{The $Q$ factor is implicitly represented by $Q_{00}$,
  $Q_{01}$, $\dots$, $Q_{0,P-1}$, and is in slow memory}
\end{algorithmic}
\end{algorithm}

Sequential TSQR on a flat tree performs the same number of flops as
sequential Householder QR, namely $2mn^2 - \frac{2n^3}{3}$ flops.
However, it performs less communication than Householder QR, as it
will be discussed in Section~\ref{pr:TSQR:other}.  Sequential TSQR
transfers $2mn - \frac{n(n+1)}{2} + \frac{mn^2}{\widetilde{W}}$ words
between slow and fast memory, in which 
\[
\widetilde{W} = W - n(n+1)/2, 
\]
and performs $\frac{2mn}{\tilde{W}}$ transfers between slow and fast
memory.  Thus the runtime for sequential TSQR is

\begin{equation}
\label{eqn:TSQR_seq_runtime}
\text{Time}_{\text{Seq.\ TSQR}}(m,n,W) = 
\left(
  2mn^2 - \frac{2n^3}{3}
\right) \gamma + 
\left( 2mn - \frac{n(n+1)}{2} + \frac{mn^2}{\widetilde{W}} \right) \beta + 
\left(
  \frac{2mn}{\widetilde{W}}
\right) \alpha 
\; \; .
\end{equation}
We note that $\widetilde{W} \stackrel{>}{\approx} 2 W / 3$, so that
the number of messages $2mn / \widetilde{W} \stackrel{<}{\approx} 3mn
/ W$.

The parallel and sequential TSQR algorithm are performed in place.
During TSQR, in the lower trapezoidal $m/P \times n$ matrix, processor
$i$ stores the Householder vectors corresponding to the local QR
factorization of its leaf node.  In the upper triangular part, it
stores first the $R_{i0}$ matrix corresponding to the local QR
factorization.  For each level $k$ of the tree at which processor $i$
participates, it will store the $R$ factor at this level.  At the last
QR factorization in which processor $i$ is involved, it will store the
Householder vectors corresponding to this QR factorization.

%% file: pr-tsqr-other.tex
\section{Other ``tall skinny'' QR algorithms}
\label{pr:TSQR:other}

There are many other algorithms besides TSQR for computing the QR
factorization of a tall and skinny matrix.  They differ in terms of
performance, flops, and accuracy, and may store the $Q$ factor in different
ways that favor certain applications over others.  In this section, we
briefly discuss the performance and summarize the numerical accuracy
of the following competitors to TSQR:
\begin{itemize}
\item variants of Gram-Schmidt
\item CholeskyQR (see \cite{stwu:02})
\item Householder QR, with a block row layout
\end{itemize}
Gram-Schmidt has two commonly used variations: ``classical'' (CGS) and
``modified'' (MGS).  Both versions have the same floating-point
operation count; MGS is more stable than CGS.
Note that we are using the row-oriented version of MGS.
CholeskyQR consists of computing the Cholesky factorization $R^T R$ of $A^T A$, and then forming $Q := A R^{-1}$.

For Householder QR, we base our parallel model on a right-looking
blocked Householder as in the ScaLAPACK
routine \lstinline!PDGEQRF!. The sequential model is based on left-looking
blocked Householder
as in the out-of-core ScaLAPACK routine
\lstinline!PFDGEQRF!~\cite{dazevedo1997design_paper}.
In the out-of core case, 
left-looking is favored instead of right-looking
in order to minimize the number of writes to slow memory (the total
amount of data moved between slow and fast memory is the same for both
left-looking and right-looking blocked Householder QR).

\begin{table}
  \centering
  \begin{tabular}{l|c|c|c}
    Parallel algorithm & \# flops & \# messages & \# words \\ \hline
    TSQR & $\frac{2mn^2}{P} + \frac{2n^3}{3} \log(P)$
         & $\log(P)$
         & $\frac{n^2}{2} \log(P)$ \\
    \lstinline!PDGEQRF!
        & $\frac{2mn^2}{P} - \frac{2n^3}{3P}$ 
        & $2n \log(P)$
        & $\frac{n^2}{2} \log(P)$ \\
    MGS by row & $\frac{2mn^2}{P}$
          & $2n \log(P)$
          & $\frac{n^2}{2}\log(P)$ \\
    CGS & $\frac{2mn^2}{P}$
          & $2n \log(P)$ 
          & $\frac{n^2}{2}\log(P)$ \\ 
    CholeskyQR & $\frac{mn^2}{P} + \frac{n^3}{3}$ 
               & $\log(P)$ 
               & $\frac{n^2}{2}\log(P)$ \\ 
  \end{tabular}
  \caption{Performance model of various parallel QR factorization
    algorithms.  Lower-order terms omitted. All parallel terms are
    counted along  the critical path. Only the  $R$ factor is computed. 
    (The $Q$ factor might be stored implicitly, explicitly or not at all depending on the algorithm.}
  \label{pr:tbl:QR:perfcomp:par}
\end{table}

\begin{table}
\small
  \centering
  \begin{tabular}{l|c|c|c}
    Sequential algorithm & \# flops & \# messages & \# words  \\\hline
    TSQR & $2mn^2 - \frac{2n^3}{3}$ 
         & $\frac{2mn}{\widetilde{W}}$
         & $2mn - \frac{n(n+1)}{2} 
           + \frac{mn^2}{\widetilde{W}}$ \\
    \lstinline!PFDGEQRF!
         & $2mn^2 - \frac{2 n^3}{3}$
         & $\frac{2mn}{W} + \frac{mn^2}{2W}$
         & $\frac{m^2 n^2}{2W} - \frac{mn^3}{6W}
            + \frac{3mn}{2} - \frac{3n^2}{4}$ \\
    MGS by row & $2mn^2$ 
        & $\frac{2mn^2}{\widetilde{W}}$
        & $\frac{3mn}{2} + \frac{m^2 n^2}{2 \widetilde{W}}$ \\
    CholeskyQR & $mn^2 + \frac{n^3}{3}$
               & $\frac{6mn}{W}$ 
               & $3mn$ \\
  \end{tabular}
  \caption{Performance model of various sequential QR factorization
    algorithms.  \lstinline!PFDGEQRF! is our model of ScaLAPACK's
    out-of-DRAM QR factorization; $W$ is the fast memory size, and
    $\tilde{W} = W - n(n+1)/2$.  Lower-order terms omitted.
    Only the $R$ factor is computed. 
    (The $Q$ factor might be stored implicitly, explicitly or not at all depending on the algorithm.}
  \label{tbl:QR:perfcomp:seq}
\end{table}

Table \ref{pr:tbl:QR:perfcomp:par} compares the performance of all the
parallel QR factorizations discussed here, and Table
\ref{tbl:QR:perfcomp:seq} compares the performance of their respective
sequential implementations, including our modeled version of
\lstinline!PFDGEQRF!.  These tables show that CholeskyQR should have
better performance than all the other methods.  This is because
CholeskyQR requires only one all-reduction operation \cite{stwu:02}.
In the parallel case, it requires $\log_2 P$ messages, where $P$ is
the number of processors.  In the sequential case, it reads the input
matrix only once.  Thus, it is optimal in the same sense that TSQR is
optimal.  Furthermore, the reduction operator is matrix-matrix
addition rather than a QR factorization of a matrix with comparable
dimensions, so CholeskyQR should always be more efficient than TSQR.
Section~\ref{pr:SS:binTSQR} supports this claim with performance data
on a cluster and a BlueGene/L platform.

However, numerical accuracy is also an important consideration for
many users.  Unlike CholeskyQR, CGS, or MGS, Householder QR is
\emph{unconditionally stable}.  That is, the computed $Q$ factors are
always orthonormal to machine precision, regardless of the properties
of the input matrix \cite{govl:96}.  This also holds for TSQR, because
the algorithm is composed entirely of no more than $P$ Householder QR
factorizations, in which $P$ is the number of input blocks.  Each of
these factorizations is itself unconditionally stable.  In contrast,
the orthogonality of the $Q$ factor computed by CGS, MGS, or
CholeskyQR depends on the condition number of the input matrix.  For
example, in CholeskyQR, the loss of orthogonality of the computed $Q$
factor depends quadratically on the condition number of the input
matrix.

However, sometimes some loss of accuracy can be tolerated, either to
improve performance, or for the algorithm to have a desirable
property.  For example, in some cases the input vectors are
sufficiently well-conditioned to allow using CholeskyQR, and the
accuracy of the orthogonalization is not so important.

We care about stability for two reasons.  First, an important
application of TSQR is the orthogonalization of basis vectors in
Krylov methods.  When using Krylov methods to compute eigenvalues of
large, ill-conditioned matrices, the whole solver can fail to converge
or have a considerably slower convergence when the orthogonality of
the Ritz vectors is poor \cite{lehoucqORTH,andrewORTH}.  Second, we
will use TSQR in Section \ref{pr:S:CAQR} as the panel factorization in
a QR decomposition algorithm for matrices of general shape.  Users who
ask for a QR factorization generally expect it to be numerically
stable.  This is because of their experience with Householder QR,
which does more work than Gaussian elimination, but produces more stable
results.  Users who are not willing to spend this additional work
already favor faster but less stable algorithms.

%% file: pr-caqr.tex
\section{Parallel CAQR}
\label{pr:S:CAQR}

The parallel CAQR (``Communication-Avoiding QR'') algorithm uses
parallel TSQR to perform a right-looking QR factorization of a dense
matrix $A$ on a two-dimensional grid of processors $P = P_r \times
P_c$.  The $m \times n$ matrix is distributed using a 2-D block cyclic
layout over the processor grid, with blocks of dimension $b \times b$.
For the sake of simplicity, we assume that all the blocks are of the
same size and square, so that they are $b \times b$; we also assume
that $m \geq n$.

In summary, the number of arithmetic
operations and words transferred is roughly the same between parallel
CAQR and ScaLAPACK's parallel QR factorization, but the number of
messages is a factor $b$ times lower for CAQR.  There is also an
analogous sequential version of CAQR, which we describe in detail in
the technical report~\cite{TSQR_technical_report}.

CAQR is based on TSQR in order to minimize communication.  At each
step of the factorization, TSQR is used to factor a panel of columns,
and the resulting Householder vectors are applied to the rest of the
matrix.  The block column QR factorization as performed in
\lstinline!PDGEQRF! is the latency bottleneck of the current ScaLAPACK
QR algorithm.  Replacing this block column factorization with TSQR,
and adapting the rest of the algorithm to work with TSQR's
representation of the panel $Q$ factors, removes the bottleneck.  We
use the reduction-to-one-processor variant of TSQR on a binary tree,
as the panel's $R$ factor need only be stored on one processor (the
processor owning the diagonal block).

CAQR is defined iteratively.  We assume that the first $j-1$
iterations of the CAQR algorithm have been performed.  That is, $j-1$
panels of width $b$ have been factored and the trailing matrix has
been updated.  The active matrix at step $j$ (that is, the part of the
matrix which needs to be worked on) is of dimension
\[
(m - (j-1) b) \times (n - (j-1) b) = m_j \times n_j.
\]

\begin{figure}[htbp]
  \begin{center}
    \mbox{
      \subfigure[Step $j$ of CAQR]{\includegraphics[scale=0.4]{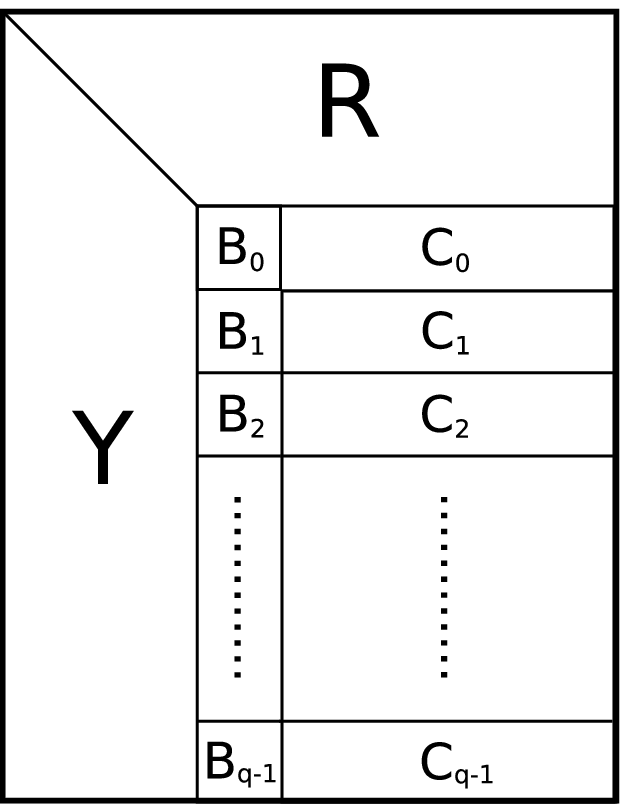}\label{fig:qr2d}} 
      \subfigure[TSQR reduction tree]{\includegraphics[scale=0.3, angle=90]{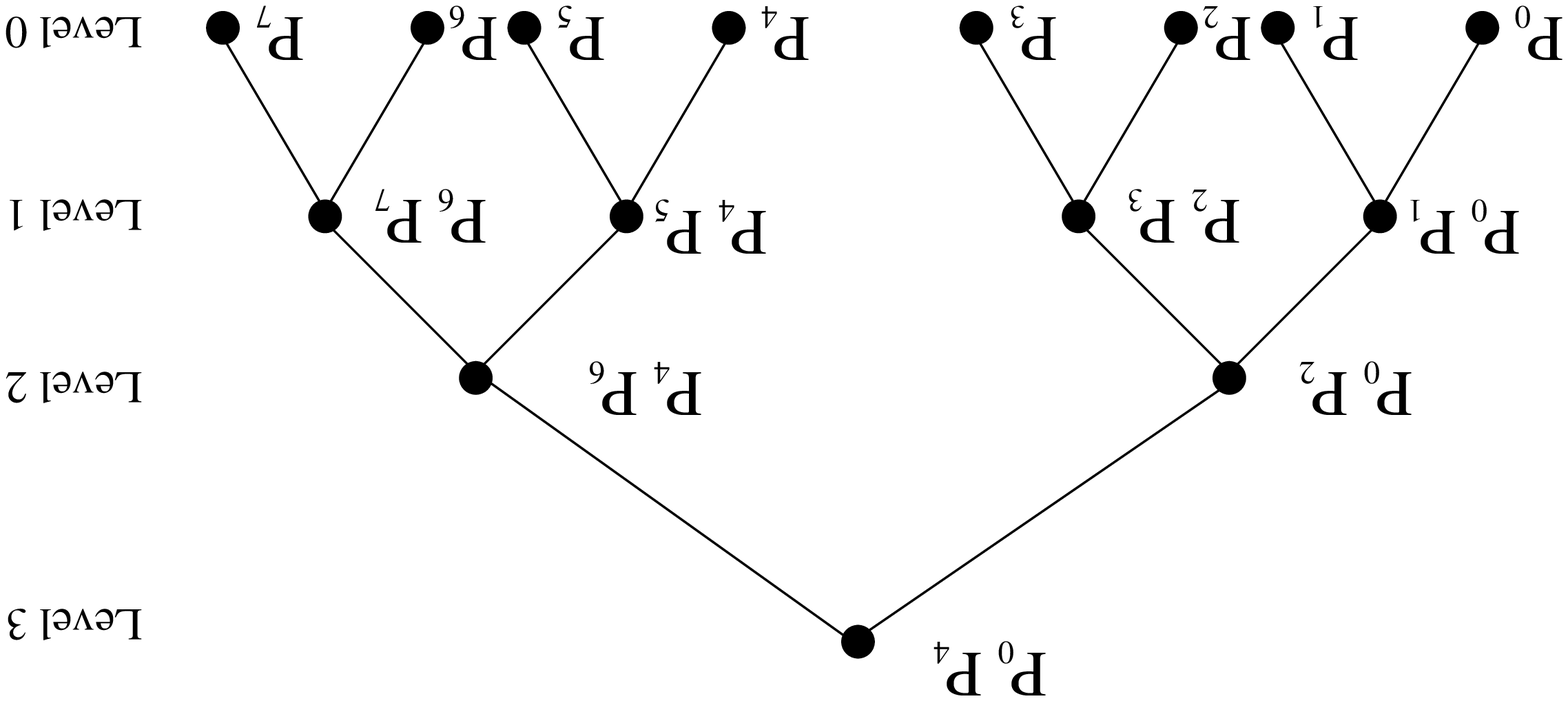}\label{fig:binTree8}}}
  \end{center}
  \caption{Step $j$ of CAQR factorization (a), and an example of a
    binary TSQR reduction tree with $8$ processors (b).  First, the
    current panel of width $b$, $B = [B_0 ; B_1 ; \cdots ; B_{q-1}]$
    is factorized using TSQR.  Here, $q$ is the number of blocks in
    the current panel.  Second, the trailing matrix, $C = [C_0 ; C_1 ;
    \cdots ; C_{q-1}]$, is updated. }
\end{figure}


Figure \ref{fig:qr2d} shows the execution of the QR factorization.
For the sake of simplicity, we suppose that processors $0$, $\dots$,
$P_r - 1$ lie in the column of processes that hold the current panel
$j$, and that $P_r$ is a power of 2.
The $m_j \times b$ matrix $B = [B_0 ; B_1 ; \cdots ; B_{q-1}]$
represents the current panel $j$.  The $m_j \times (n_j - b)$ matrix
$C = [C_0 ; C_1 ; \cdots ; C_{q-1}]$ is the trailing matrix that needs
to be updated after the TSQR factorization of $B$.  For each processor
$i$, the first $b$ rows of its first block row of $B$ and $C$ are
$B_i$ and $C_i$ respectively.

We first introduce some notation to help us refer to different parts
of a binary TSQR reduction tree.  TSQR takes place in ($\log_2{P_r} +
1$) steps, starting from the bottom level $k = 0$ of a binary tree.
Each node of the binary tree is associated with a set of processors.
We use the following notations:

\begin{itemize}
\item $level(i,k) = \left\lfloor \frac{i}{2^k} \right\rfloor$ denotes
  the node at level $k$ of the reduction tree which is assigned to a
  set of processors that includes processor $i$.  
\item $first\_proc(i,k) = 2^k level(i,k)$ is the index of the
  ``first'' processor associated with the node $level(i,k)$ at stage
  $k$ of the reduction tree.  In a reduction (not an all-reduction),
  it receives the messages from its neighbors and performs the local
  computation.
\item $target\_first\_proc(i,k) = first\_proc(i,k) + 2^{k-1}$ is the
  index of the processor with which $first\_proc(i,k)$ exchanges data
  in a reduction at level $k$.
\end{itemize}

A binary TSQR reduction tree for $8$ processors is shown in
Figure~\ref{fig:binTree8}. For example, the processors $P_4$ and $P_6$
are affected to the right node at level $k=2$.  With the above
notation, the processors in the range $i=4, \ldots, 7$ can compute
easily the two processors affected to this node, that is
$first\_proc(i,2) = 4$ and $target\_first\_proc(i,2) = 6$.

Algorithm \ref{Alg:CAQR:j} outlines the right-looking parallel QR
decomposition.  At iteration $j$, first, the block column $j$ is
factored using TSQR.  After the block column factorization is
complete, the trailing matrix is updated as follows.  The update
corresponding to the QR factorization at the leaves of the TSQR tree
is performed locally on every processor.  The updates corresponding to
the upper levels of the TSQR tree are performed between groups of
neighboring trailing matrix processors.  Note that only one of the
trailing matrix processors in each neighbor group continues to be
involved in successive trailing matrix updates.  This allows overlap
of computation and communication, as the uninvolved processors can
finish their computations in parallel with successive reduction
stages.

\begin{algorithm}[h!]
\caption{Right-looking parallel CAQR factorization}
\label{Alg:CAQR:j}
\begin{algorithmic}[1]
\For{$j = 1$ to $n/b$}
  \State{The column of processors that holds panel $j$ computes a TSQR
    factorization of this panel.  The Householder vectors are stored
    in a tree-like structure as described in Section 
    \ref{S:TSQR:impl}.}\label{Alg:CAQR:j:local-factor}
  \State{Each processor $p$ that belongs to the column of processes
    holding panel $j$ broadcasts along its row of processors the $m_j
    / P_r \times b$ rectangular matrix that holds the two sets of
    Householder vectors.  Processor $p$ also broadcasts two arrays of
    size $b$ each, containing the Householder multipliers $\tau_p$.}
  \State{Each processor in the same process row as processor $p$, $0
    \leq p < P_r$, forms $T_{p0}$ and updates its local trailing
    matrix $C$ using $T_{p0}$ and $Y_{p0}$. (This computation involves
    all processors.)}\label{Alg:CAQR:j:local-update}
  \For{$k = 1$ to $\log P_r$, the processors that lie in the same row 
    as processor $p$, where $0 \leq p < P_r$ equals 
    $first\_proc(p,k)$ or $target\_first\_proc(p,k)$,
    respectively.}
      \State{Processors in the same process row as
        $target\_first\_proc(p,k)$ form $T_{level(p,k),k}$ locally.  
	They also compute local pieces of 
	$W = Y_{level(p,k),k}^T C_{target\_first\_proc(p,k)}$,
        leaving the results distributed.  This computation is 
        overlapped with the communication in Line \ref{step_comm1}.}\label{Alg:CAQR:j:overlap1}

      \State{Each processor in the same process row as 
        $first\_proc(p,k)$ sends to the processor in the 
        same column and belonging to the row of processors of 
        $target\_first\_proc(p,k)$ the local pieces of
        $C_{first\_proc(p,k)}$.}\label{step_comm1}

      \State{Processors in the same process row as 
        $target\_first\_proc(p,k)$ compute local pieces of
        \[
        W = T_{level(p,k),k}^T \left( C_{first\_proc(p,k)} + W \right).
        \]}

      \State{Each processor in the same process row as 
        $target\_first\_proc(p,k)$ sends to the processor 
        in the same column and belonging to the process row
        of $first\_proc(p,k)$ the local pieces of $W$.}\label{step_comm2}

      \State{Processors in the same process row as
        $first\_proc(p,k)$ and
        $target\_first\_proc(p,k)$ each complete the rank-$b$
        updates $C_{first\_proc(p,k)} := C_{first\_proc(p,k)} - W$ and 
        $C_{target\_first\_proc(p,k)} := C_{target\_first\_proc(p,k)} -
        Y_{level(p,k),k} \cdot W$ locally.  The latter computation
        is overlapped with the communication in Line \ref{step_comm2}.}\label{Alg:CAQR:j:overlap2}
  \EndFor
\EndFor
\end{algorithmic}
\end{algorithm}

We see that CAQR consists of $\frac{n}{b}$ TSQR factorizations
involving $P_r$ processors each, and $n/b - 1$ applications of the
resulting Householder vectors.  Table \ref{pr:tbl:CAQR:par:model}
expresses the performance model over a rectangular $P_r \times P_c$
grid of processors.  A detailed derivation of the model is given
in~\cite{TSQR_technical_report}.  According to the table, the number
of arithmetic operations and words transferred is roughly the same
between parallel CAQR and ScaLAPACK's parallel QR factorization, but
the number of messages is a factor $b$ times lower for CAQR.

\begin{table}[h]
\small
\centering
\begin{tabular}{l | l}
            & Parallel CAQR \\ \hline
\# messages & $\frac{3n}{b} \log P_r + \frac{2n}{b} \log P_c$ \\ \hline
\# words    & $\left( 
                   \frac{n^2}{P_c} 
                   + \frac{bn}{2} 
               \right) \log P_r
               + \left( 
                   \frac{mn - n^2/2}{P_r} + 2n 
               \right) \log P_c$ \\ \hline
\# flops    & $\frac{2n^2(3m-n)}{3P} 
               + \frac{bn^2}{2P_c} 
               + \frac{3bn(2m - n)}{2P_r} 
               + \left( \frac{4 b^2 n}{3} 
                   + \frac{n^2 (3b+5)}{2 P_c} 
               \right) \log P_r
               - b^2 n$ \\ \hline
             & ScaLAPACK's \texttt{PDGEQRF} \\ \hline
\# messages & $3n \log P_r + \frac{2n}{b} \log P_c$ \\ \hline
\# words    & $\left( 
                   \frac{n^2}{P_c} 
                   + bn 
               \right) \log P_r
               + \left( 
                   \frac{mn - n^2/2}{P_r} 
                   + \frac{bn}{2} 
               \right) \log P_c$ \\ \hline
\# flops    & $\frac{2n^2(3m-n)}{3P} 
               + \frac{bn^2}{2P_c} 
               + \frac{3bn(2m - n)}{2P_r}
               - \frac{b^2 n}{3 P_r}$
            \\ \hline
\end{tabular}
\caption{Performance models of parallel CAQR and ScaLAPACK's
  \lstinline!PDGEQRF! when factoring an $m \times n$ matrix,
  distributed in a 2-D block cyclic layout on a $P_r \times P_c$ grid
  of processors with square $b \times b$ blocks.  All terms are
  counted along the critical path.  In this table, ``flops'' only
  includes floating-point additions and multiplications, not
  floating-point divisions.  Some lower-order terms are omitted.  We
  generally assume $m \geq n$.}
\label{pr:tbl:CAQR:par:model}
\end{table}

The parallelization of the computation is represented by the number of
flops in Table \ref{pr:tbl:CAQR:par:model}.  The first, dominant, 
term for CAQR
represents mainly the parallelization of the local Householder update
corresponding to the leaves of the TSQR tree (the matrix-matrix
multiplication in line \ref{Alg:CAQR:j:local-update} of Algorithm
\ref{Alg:CAQR:j}), and matches the first term for \lstinline!PDGEQRF!.
The second term for CAQR corresponds to forming the $T_{p0}$ matrices
for the local Householder update in line \ref{Alg:CAQR:j:local-update}
of the algorithm, and also has a matching term for
\lstinline!PDGEQRF!.  The third term for CAQR represents the QR
factorization of a panel of width $b$ that corresponds to the leaves
of the TSQR tree (part of line \ref{Alg:CAQR:j:local-factor}) and part
of the local rank-$b$ update (triangular matrix-matrix multiplication)
in line \ref{Alg:CAQR:j:local-update} of the algorithm, and also has a
matching term for \lstinline!PDGEQRF!.

The fourth, lower order, term in the number of flops for CAQR represents the
redundant computation introduced by the TSQR formulation. In this
term, the number of flops performed for computing the QR factorization
of two upper triangular matrices at each node of the TSQR tree is
$(2/3) nb^2 \log(P_r)$.  The number of flops performed during the
Householder updates issued by each QR factorization of two upper
triangular matrices is $n^2 (3b+5)/(2 P_c) \log(P_r)$.

We note that standard optimizations like overlapping computation and
communication, as in look-ahead techniques, are possible with CAQR.
With the look-ahead right-looking approach, the communications are
pipelined from left to right.  At each step of factorization, we would
model the latency cost of the broadcast within rows of processors as
$2$ instead of $\log P_c$.  Also, the runtime estimation in Table
\ref{pr:tbl:CAQR:par:model} does not take into account the overlap of
computation and communication in lines \ref{Alg:CAQR:j:overlap1} and
\ref{step_comm1} or in lines \ref{step_comm2} and
\ref{Alg:CAQR:j:overlap2} of Algorithm \ref{Alg:CAQR:j}.  Suppose that
at each step of the QR factorization, the condition
\[
     \alpha + \beta \frac{b (n_j - b)}{P_c} 
     > 
     \gamma b (b + 1) \frac{n_j - b}{P_c}
\]
is fulfilled,  this is the case for example when $\beta / \gamma >
b+1$,  then the fourth flops term that accounts for the redundant
computation is decreased by $n^2 (b+1) \log(P_r) / P_c$, about a
factor of $3$.

%% file: pr-experiments.tex
\section{Experimental results}
\label{pr:experiments}

In this section we present the performance of sequential and parallel
TSQR on several computational systems.  We also use the performance
model of CAQR in Table~\ref{pr:tbl:CAQR:par:model} to predict its
performance and compare it to PDGEQRF.  The actual implementation and
measurements of parallel CAQR are currently underway.

TSQR (and its associated CAQR factorization algorithm on a 2-D matrix
layout) is not a single algorithm, but a space of possible algorithms.
It encompasses all possible reduction tree shapes, including:
\begin{enumerate}
  \item Binary (to minimize number of messages in the parallel case)
  \item Flat (to minimize communication volume in the sequential case)
  \item Hybrid (to account for network topology, and/or to balance
    bandwidth demands with maximum parallelism)
\end{enumerate}
as well as all possible ways to perform the local QR factorizations,
including:
\begin{enumerate}
\item (Possibly multithreaded) standard LAPACK (DGEQRF)
  \item An existing parallel QR factorization, such as ScaLAPACK's
    PDGEQRF 
  \item A recursive QR factorization (e.g., \cite{elmroth1998new,elmroth2000applying})
\end{enumerate}
Choosing the right combination of parameters can help minimize
communication between any or all of the levels of the memory
hierarchy, from cache and shared-memory bus, to DRAM and local disk,
to parallel filesystem and distributed-memory network interconnects,
to wide-area networks.

The huge tuning space makes it a challenge to pick the right platforms
for experiments.  Luckily, TSQR's hierarchical structure makes tunings
\emph{composable}.  For example, once we have a good choice of
parameters for TSQR on a single multicore node, we don't need to
change them when we tune TSQR for a cluster of these nodes.  From the
cluster perspective, it is as if the performance of the individual
nodes improved.  This means that we can benchmark TSQR on a small,
carefully chosen set of scenarios, with confidence that they represent
many platforms of interest.

Previous work covers some parts of the tuning space.  Gunter et al.\
implement an out-of-DRAM version of TSQR on a flat tree, and use a
parallel distributed QR factorization routine to factor in-DRAM blocks
\cite{gunter2005parallel}.  Buttari et al.\ suggest using a QR
factorization of this type to improve performance of parallel QR on
commodity multicore processors \cite{buttari2007class}.  Quintana-Ort\'i
et al.\ develop two variations on block QR factorization algorithms,
and use them with a dynamic task scheduling system to parallelize the
QR factorization on shared-memory machines
\cite{quintana-orti2008scheduling}.  Kurzak and Dongarra use similar
algorithms, but with static task scheduling, to parallelize the QR
factorization on Cell processors \cite{kurzak2008qr}.  Pothen and
Raghavan \cite{pothen1989distributed} and Cunha et al.
\cite{cunha2002new} both benchmarked parallel TSQR using a binary tree
on a distributed-memory cluster, and implemented the local QR
factorizations with a single-threaded version of \lstinline!DGEQRF!.
All these researchers observed significant performance improvements
over previous QR factorization algorithms. 
The only parallel implementations of CAQR we are aware of are parallel CAQR
with a flat tree in the shared memory context.  These have recently been
presented in~\cite{buttari2007class,quintana-orti2008scheduling}.  To our
knowledge, there is no implementation of parallel CAQR in the distributed
context (neither flat tree nor binary tree).

We choose to run two sets of experiments for TSQR.  The first set
covers the out-of-DRAM case on a single CPU, and the results are
presented in Section~\ref{pr:SS:flatTSQR}.  We use a laptop with a
single PowerPC CPU for these experiments.  The second set, presented
in Section~\ref{pr:SS:binTSQR}, is like the parallel experiments of
previous authors in that it uses a binary tree on a distributed-memory
cluster, but it improves on their approach by using a better local QR
factorization (the recursive approach -- see
\cite{elmroth1998new,elmroth2000applying}).  We use two distributed-memory machines:
a Pentium III cluster (``Beowulf'') and a BlueGene/L (``Frost'').

In Section \ref{pr:S:CAQR:perfest}, we estimate performance of
parallel CAQR on our projection of a future petascale machine with
8192 processors (``Peta'').  Detailed performance evaluation on two
different parallel machines, an existing IBM POWER5 and a grid formed
by 128 processors linked together by the Internet, can be found in the
technical report~\cite{TSQR_technical_report}.

\subsection{Tests of sequential TSQR on a flat tree}
\label{pr:SS:flatTSQR}
We developed an out-of-DRAM version of TSQR that uses a flat reduction
tree.  It invokes the system vendor's native BLAS and LAPACK
libraries.  Thus, it can exploit a multithreaded BLAS on a machine
with multiple CPUs, but the parallelism is limited to operations on a
single block of the matrix.  We used standard POSIX blocking file
operations, and made no attempt to overlap communication and
computation.  Exploiting overlap could at best double the performance.

We ran sequential tests on a laptop with a single PowerPC CPU.
Details of the platform are as follows:
\begin{itemize}
\item Single-core PowerPC G4 (1.5 GHz), with 512 KB of L2 cache, 512
  MB of DRAM on a 167 MHz bus, One Fujitsu MHT2080AH 80 HB hard drive
  (5400 RPM).
\end{itemize}

In our experiments, we first used both out-of-DRAM TSQR and standard
LAPACK QR to factor a collection of matrices that use only slightly
more than half of the total DRAM for the factorization.  This was so
that we could collect comparison timings.  Then, we ran only
out-of-DRAM TSQR on matrices too large to fit in DRAM or swap space,
so that an out-of-DRAM algorithm is necessary to solve the problem at
all.  For the latter timings, we extrapolated the standard LAPACK QR
timings up to the larger problem sizes, in order to estimate the
runtime if memory were unbounded.  LAPACK's QR factorization swaps so
much for out-of-DRAM problem sizes that its actual runtimes are many
times larger than these extrapolated unbounded-memory runtime
estimates. Note that once an in-DRAM algorithm begins swapping, it
becomes so much slower that most users prefer to abort the computation
and try solving a smaller problem.

We used the following power law for the extrapolation:
\[
t = A_1 b m^{A_2} n^{A_3},
\]
in which $t$ is the time spent in computation, $b$ is the number of
input matrix blocks, $m$ is the number of rows per block, and $n$ is
the number of columns in the matrix.  After taking logarithms of both
sides, we performed a least squares fit of $\log(A_1)$, $A_2$, and
$A_3$.  The value of $A_2$ was 1, as expected.  The value of $A_3$ was
about 1.6.  This is less than 2 as expected, given that increasing the
number of columns increases the computational intensity and thus the
potential for exploitation of locality (a Level 3 BLAS effect).  We expect
around two digits of accuracy in the parameters, which in themselves
are not as interesting as the extrapolated runtimes; the parameter
values mainly serve as a sanity check.

\begin{figure}
  \begin{center}
    \mbox{
      \subfigure[Measured data]{\includegraphics[scale=0.3]{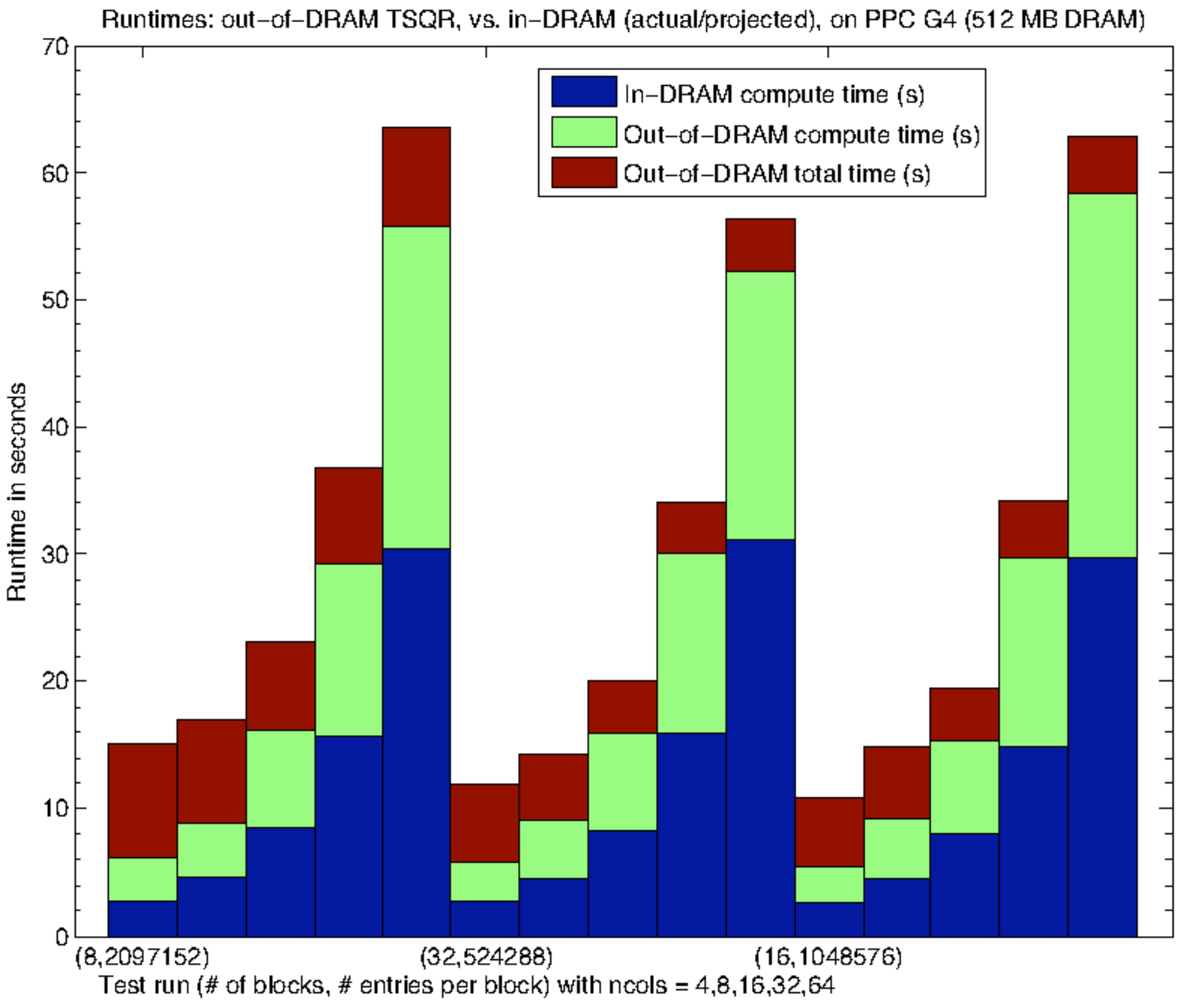}\label{fig:tsqr-laptop-no-extrap}} 
      \subfigure[Extrapolated runtime]{\includegraphics[scale=0.3]{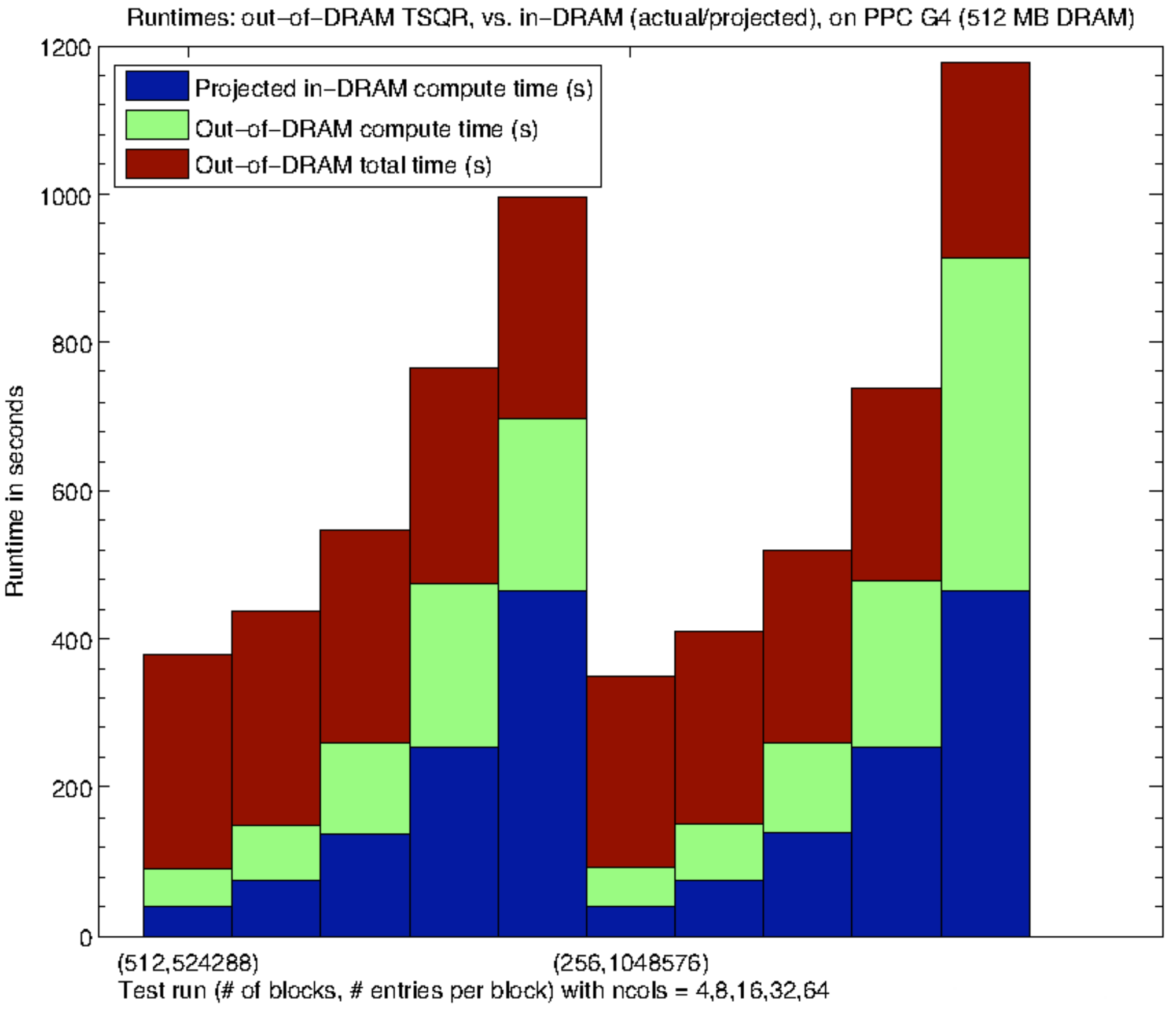}\label{fig:tsqr-laptop-extrap}}}
  \end{center}
  \caption{\small Runtimes (in seconds) of out-of-DRAM TSQR compared
    against (a) measured data and (b) extrapolated runtime of standard
    QR (LAPACK's \lstinline!DGEQRF!) on a single-processor laptop.
    For the measured data, we limit memory usage to 256 MB, which is
    half of the laptop's total system memory, so that we can collect
    performance data for DGEQRF. For extrapolated data, we use the
    measured data to construct a power-law performance extrapolation.
    The graphs show different choices of block dimensions and number
    of blocks $P$.  The top of the blue bar is (a) the benchmarked total
    runtime for \lstinline!DGEQRF! and (b) the extrapolated total
    runtime for \lstinline!DGEQRF!, the top of the green bar is the
    benchmarked compute time for TSQR, and the top of the brown bar is
    the benchmarked total time for TSQR.  Thus the height of the brown
    bar alone is the I/O time.  Note that LAPACK starts and ends in
    DRAM (if it could fit in DRAM), and TSQR starts and ends on disk.}
\end{figure}

Figure~\ref{fig:tsqr-laptop-no-extrap} shows the measured in-DRAM
results on the laptop platform, and
Figure~\ref{fig:tsqr-laptop-extrap} shows the (measured TSQR,
extrapolated LAPACK) out-of-DRAM results on the same platform.  In
these figures, the amount of memory, and so the total number of
matrix entries is constant for all the experiments: $m \cdot n = 2^24$.
This means the total volume of communication is the same for all
experiments.
The number of blocks $P$ used, and so the number of
matrix entries  per block $mn/P$, is the same for each group of five bars,
and is shown in a label under the horizontal axis.
Within each group of 5 bars, we varied the number of matrix columns to be
4, 8, 16, 32, and 64.
Note that we
have not tried to overlap I/O and computation in this implementation.
The trends in Figure \ref{fig:tsqr-laptop-no-extrap} suggest that the
extrapolation is reasonable: TSQR takes about twice as much time for
computation as does standard LAPACK QR, and the fraction of time spent
in I/O is reasonable and decreases with problem size.

TSQR assumes that the matrix starts and ends on disk, whereas LAPACK
starts and ends in DRAM.  Thus, to compare the two, one could also
estimate LAPACK performance with infinite DRAM but where the data
starts and ends on disk.  The height of the reddish-brown bars in
Figures \ref{fig:tsqr-laptop-no-extrap} and
\ref{fig:tsqr-laptop-extrap} is the I/O time for TSQR, which can be
used to estimate the LAPACK I/O time.  
This is reasonable since the volume of communication in the two
cases is the same, and the fact that the reddish-brown bars are
of similar height for different values of $P$, 
shows that the communication is bandwidth dominated.
Add this to the blue bar (the
LAPACK compute time) to estimate the runtime when the LAPACK QR
routine must load the matrix from disk and store the results back to
disk.

The main purpose of our out-of-DRAM code is not to outperform existing
in-DRAM algorithms, but to be able to solve classes of problems which
the existing algorithms cannot solve.  The above graphs show that the
penalty of an explicitly swapping approach is about 2x, which is small
enough to warrant its practical use.  This holds even for problems
with a relatively low computational intensity, such as when the input
matrix has very few columns.  Furthermore, picking the number of
columns sufficiently large may allow complete overlap of file I/O by
computation.

\subsection{Tests of parallel TSQR on a binary tree}
\label{pr:SS:binTSQR}
\begin{table}
  \begin{center}
{\small
\begin{tabular}{c|r|r|r|r|r|r} 
\# procs & CholeskyQR & TSQR & CGS & MGS & TSQR & ScaLAPACK \\
 & & (\lstinline!DGEQR3!) & & & (\lstinline!DGEQRF!) & (\lstinline!PDGEQRF!) \\ \hline
  1 & 1.02 & 4.14 &  3.73 &  7.17 &  9.68 & 12.63 \\ 
  2 & 0.99 & 4.00 &  6.41 & 12.56 & 15.71 & 19.88 \\ 
  4 & 0.92 & 3.35 &  6.62 & 12.78 & 16.07 & 19.59 \\ 
  8 & 0.92 & 2.86 &  6.87 & 12.89 & 11.41 & 17.85 \\ 
 16 & 1.00 & 2.56 &  7.48 & 13.02 &  9.75 & 17.29 \\ 
 32 & 1.32 & 2.82 &  8.37 & 13.84 &  8.15 & 16.95 \\ 
 64 & 1.88 & 5.96 & 15.46 & 13.84 &  9.46 & 17.74 \\ 
\end{tabular}
}
\caption{Runtime in seconds of various parallel QR factorizations on
  the Beowulf machine.  The total number of rows $m = 100000$ and the
  ratio $\lceil n / \sqrt{P} \rceil = 50$ (with $P$ being the number
  of processors) were kept constant as $P$ varied from 1 to 64.  This
  illustrates weak scaling with respect to the square of the number of
  columns $n$ in the matrix, which is of interest because the number
  of floating-point operations in sequential QR is $\mathcal{O}(mn^2)$.  If
  an algorithm scales perfectly, then all the runtimes in that
  algorithm's column should be constant.  Both the $Q$ and $R$ factors
  were computed explicitly; in particular, for those codes which form
  an implicit representation of $Q$, the conversion to an explicit
  representation was included in the runtime measurement.}
\label{tbl:TSQR:cluster:weak:n}
\end{center}
\end{table}

\begin{table}
  \begin{center}
{\small
\begin{tabular}{c|r|r|r|r|r|r}
\# procs & CholeskyQR & TSQR & CGS & MGS & TSQR & ScaLAPACK \\
 & & (\lstinline!DGEQR3!) & & & (\lstinline!DGEQRF!)  & (\lstinline!PDGEQRF!) \\ \hline
  1 & 0.45 & 3.43 & 3.61 &  7.13 &  7.07 &  7.26 \\ 
  2 & 0.47 & 4.02 & 7.11 & 14.04 & 11.59 & 13.95 \\ 
  4 & 0.47 & 4.29 & 6.09 & 12.09 & 13.94 & 13.74 \\ 
  8 & 0.50 & 4.30 & 7.53 & 15.06 & 14.21 & 14.05 \\ 
 16 & 0.54 & 4.33 & 7.79 & 15.04 & 14.66 & 14.94 \\ 
 32 & 0.52 & 4.42 & 7.85 & 15.38 & 14.95 & 15.01 \\ 
 64 & 0.65 & 4.45 & 7.96 & 15.46 & 14.66 & 15.33 \\ 
\end{tabular}
}
\caption{Runtime in seconds of various parallel QR factorizations on
  the Beowulf machine, illustrating weak scaling with respect to the
  total number of rows $m$ in the matrix.  The ratio $\lceil m/P
  \rceil = 100000$ and the total number of columns $n = 50$ were kept
  constant as the number of processors $P$ varied from 1 to 64.  If an
  algorithm scales perfectly, then all the runtimes in that
  algorithm's column should be constant.  For those algorithms which
  compute an implicit representation of the $Q$ factor, that
  representation was left implicit.}
\label{tbl:TSQR:cluster:weak:m}
  \end{center}
\end{table}

\begin{table}
  \begin{center}
{\small
\begin{tabular}{c|r|r|r|r|r|r}
\# procs & CholeskyQR & TSQR & CGS & MGS & TSQR & ScaLAPACK \\
 & & (\lstinline!DGEQR3!) & & & (\lstinline!DGEQRF!)  & (\lstinline!PDGEQRF!) \\ \hline
 32& 0.140  &  0.453  &  0.836  &  0.694  &  1.132  &  1.817 \\
 64& 0.075  &  0.235  &  0.411  &  0.341  &  0.570  &  0.908 \\
128& 0.038  &  0.118  &  0.180  &  0.144  &  0.247  &  0.399 \\
256& 0.020  &  0.064  &  0.086  &  0.069  &  0.121  &  0.212 \\
\end{tabular}
}
\caption{Runtime in seconds of various parallel QR factorizations on
  the Frost machine on a $10^6 \times 50$ matrix.  This metric
  illustrates strong scaling (constant problem size, but number of
  processors increases).}
\label{tbl:TSQR:BGL:strong}
  \end{center}
\end{table}

We also present results for a parallel MPI implementation of TSQR on a
binary tree.  Rather than LAPACK's \lstinline!DGEQRF!, the code uses a
custom local QR factorization, \lstinline!DGEQR3!, based on the
recursive approach of Elmroth and Gustavson
\cite{elmroth2000applying}.  Tests show that \lstinline!DGEQR3!
consistently outperforms LAPACK's \lstinline!DGEQRF! by a large margin
for matrix dimensions of interest.

We ran parallel TSQR on the following distributed-memory machines:
\begin{itemize}
\item Pentium III cluster (``Beowulf''), operated by the University of
  Colorado Denver.  It has 35
  dual-socket 900 MHz Pentium III nodes with Dolphin interconnect.
  Peak floating-point rate is 900 Mflop/s per processor, network
  latency is less than 2.7 $\mu$s, benchmarked\footnote{See
  \url{http://www.dolphinics.com/products/benchmarks.html}.}, and
  network bandwidth is 350 MB/s, benchmarked upper bound.
\item IBM BlueGene/L (``Frost''), operated by the National Center for
  Atmospheric Research.  We use one BlueGene/L rack with 1024 700 MHz
  compute CPUs.  Peak floating-point rate is 2.8 Gflop/s per
  processor, network\footnote{The BlueGene/L has two separate networks
  -- a torus for nearest-neighbor communication and a tree for
  collectives.  The latency and bandwidth figures here are for the
  collectives network.} latency is 1.5 $\mu$s, hardware, and network
  one-way bandwidth is 350 MB/s, hardware.
\end{itemize}

The experiments compare many different implementations of a parallel
QR factorization.
TSQR was tested both with the recursive local QR factorization
\lstinline!DGEQR3!, and the standard LAPACK routine
\lstinline!DGEQRF!.  Both CGS and MGS (by row) were timed.

Tables \ref{tbl:TSQR:cluster:weak:n} and \ref{tbl:TSQR:cluster:weak:m} show the
results of two different performance experiments on the Pentium III cluster.
In the first of these, the total number of rows $m = 100,000$ and the ratio
$\lceil n / \sqrt{P} \rceil = 50$ (with $P$ being the number of processors)
were kept constant as $P$ varied from $1$ to $64$.  This was meant to illustrate
weak scaling with respect to $n^2$ (the square of the number of columns in the
matrix), which is of interest because the number of floating-point operations
in sequential QR is $\mathcal{O}(mn^2)$.  If an algorithm scales perfectly, then all
the runtimes shown in that algorithm's column should be constant.  Both the
tall and skinny $Q$ and the square $R$ factors were computed explicitly; in
particular, for those codes which form an implicit representation of $Q$, the
conversion to an explicit representation was included in the runtime
measurement.  The results show that TSQR scales better than CGS or MGS (by
row), and significantly outperforms ScaLAPACK's QR.  Also, using the recursive
local QR in TSQR, rather than LAPACK's QR, more than doubles performance.
CholeskyQR gets the best performance of all the algorithms, but at the expense
of significant loss of orthogonality when the initial matrix $A$ is
ill-conditioned.  Note that, in this case ($Q$ and $R$ requested), CholeskyQR,
CGS, and MGS perform half the flops of the Householder based algorithms, TSQR-DGEQR3, TSQR-DGEQRF, and
PDGEQRF ($2mn^2$ versus $4mn^2$). 

Table \ref{tbl:TSQR:cluster:weak:m} shows the results of the second
set of experiments on the Pentium III cluster.  In these experiments,
we also illustrate weak scaling with respect to the total number of
rows $m$ in the matrix.  For this, the ratio $\lceil m/P \rceil =
100,000$ and the total number of columns $n = 50$ were kept constant as
the number of processors $P$ varied from 1 to 64.  Unlike in the
previous set of experiments, for those algorithms which compute an
implicit representation of the $Q$ factor, that representation was
left implicit.  The results show that TSQR scales well.  In
particular, when using TSQR with the recursive local QR factorization,
there is almost no performance penalty for moving from one processor
to two, unlike with CGS, MGS, and ScaLAPACK's QR.  Again, the
recursive local QR significantly improves TSQR performance; here it is
the main factor in making TSQR perform better than ScaLAPACK's QR.
Note that, in this case (only $R$ requested), CholeskyQR,
performs half the flops of all the others algorithm CGS, MGS, TSQR-DGEQR3, TSQR-DGEQRF, and
PDGEQRF ($mn^2$ versus $2mn^2$). 

Table \ref{tbl:TSQR:BGL:strong} shows the results of the third set of
experiments, which was performed on the BlueGene/L cluster ``Frost.''
These data show performance per processor (Mflop / s / (number of
processors)) on a matrix of constant dimensions $10^6 \times 50$, as
the number of processors was increased.  This illustrates strong scaling.  If
an algorithm scales perfectly, then all the numbers in that algorithm's 
column should decrease proportionally to $P$, i.e. halve from row to row.
For ScaLAPACK's QR factorization, we used
\lstinline!PDGEQRF!.
We observe that using the recursive local QR factorization with TSQR makes it
 clearly outperfom ScaLAPACK.
Note that, in this case ($Q$ and $R$ requested), CholeskyQR,
CGS, and MGS perform half the flops of the Householder based algorithms, TSQR-DGEQR3, TSQR-DGEQRF, and
PDGEQRF ($2mn^2$ versus $4mn^2$). 

Both the Pentium III and BlueGene/L platforms have relatively slow
processors with a relatively low-latency interconnect.  TSQR was
optimized for the opposite case of fast processors and expensive
communication.  Nevertheless, TSQR outperforms ScaLAPACK's QR by over
$6.7\times$ on 16 processors (and $3.5\times$ on 64 processors) on the
Pentium III cluster, and 
$4.0\times$ on 32 processors (and $3.3\times$ on 256 processors) on the
BlueGene/L machine.

\subsection{Performance estimation of parallel CAQR }
\label{pr:S:CAQR:perfest}

We use the performance model developed in Section~\ref{pr:S:CAQR} to
estimate the performance of parallel CAQR on a model of a petascale
machine.  We expect CAQR to outperform ScaLAPACK, in part because it
uses a faster algorithm for performing most of the computation of each
panel factorization (\lstinline!DGEQR3! vs.\ \lstinline!DGEQRF!), and
in part because it reduces the latency cost.  Our performance model
uses the same time per floating-point operation for both CAQR and
\lstinline!PDGEQRF!.  Hence our model evaluates the improvement due
only to reducing the latency cost.

Our projection of a future petascale machine (``Peta'') has 8192
processors.  Each ``processor'' of Peta may itself be a parallel
multicore node, but we consider it as a single fast sequential
processor for the sake of our model.  Here are the parameters we use:
\begin{itemize}
\item Peta. Peak floating-point rate is 500 Gflop/s per processor,
  network latency is 10 $\mu$s, and network bandwidth is 4 GB/s.
\end{itemize}

We evaluate the performance using matrices of size $n \times n$,
distributed over a $P_r \times P_c$ grid of $P$ processors using a 2D
block cyclic distribution, with square blocks of size $b \times b$.
We estimate the best performance of CAQR and \lstinline!PDGEQRF! for a
given problem size $n$ and a given number of processors $P$, by
finding the optimal values for the block size $b$ and the shape of the
grid $P_r \times P_c$ in the allowed ranges.  The matrix size $n$ is
varied in the range $10^3$, $10^{3.5}$, $10^4$, $\dots$, $10^{7.5}$.
The block size $b$ is varied in the range $1$, $5$, $10$, $\dots$,
$50$, $60$, $\dots$, $\min(200, m/P_r, n/P_c)$.  The number of
processors is varied from $1$ to the largest power of $2$ smaller than
$p_{max}$, in which $p_{max}$ is the maximum number of processors
available in the system.  The values for $P_r$ and $P_c$ are also
chosen to be powers of two.

When we evaluate the model, we set the floating-point performance
value in the model so that the modeled floating-point rate is 80\% of
the machine's peak rate, so as to capture realistic performance on the
local QR factorizations.  The inverse network bandwidth $\beta$ has
units of seconds per word.  The white regions in the plots signify
that the problem needed more memory than available on the machine.

\begin{figure}
  \begin{center}
    \mbox{ \subfigure[Speedup
      CAQR]{\includegraphics[scale=0.35]{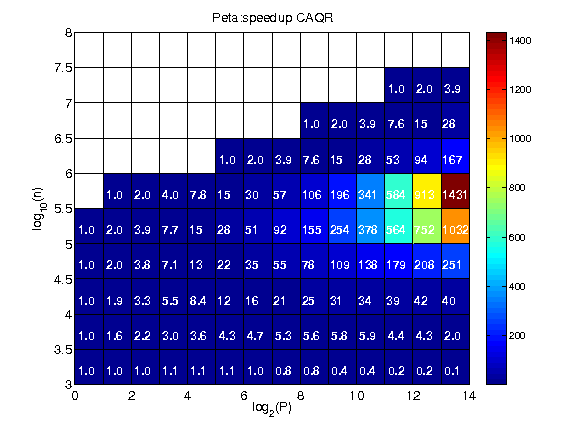}\label{peta_spdCAQR}}
      \subfigure[Speedup
      PDGEQRF]{\includegraphics[scale=0.35]{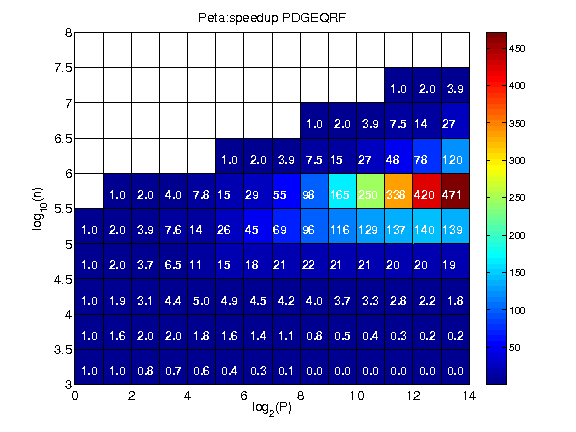}\label{peta_spdQRF}}
      } \subfigure[Comparison]
      {\includegraphics[scale=0.35]{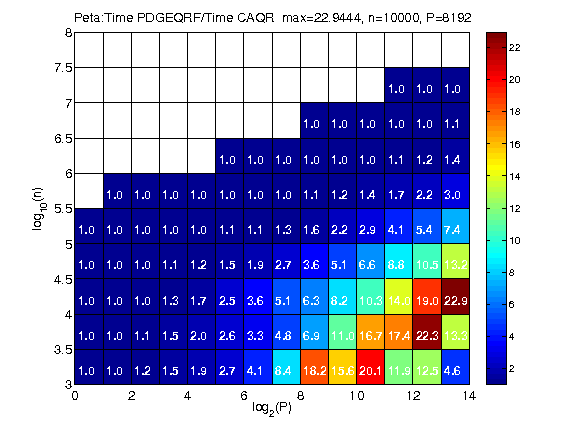}\label{peta_cmp}}
  \end{center}
  \caption{\label{fig:PerfComp_Peta}Performance prediction comparing
    CAQR and PDGEQRF on Peta.}
  \label{peta_all}
\end{figure}

Figure~\ref{fig:PerfComp_Peta} shows our performance estimates of CAQR
and \lstinline!PDGEQRF! on the Petascale machine, in which we display
\begin{itemize}
\item Figure~\ref{peta_spdCAQR} -- the best speedup obtained by CAQR,
  with respect to the runtime using the fewest number of processors
  with enough memory to hold the matrix (which may be more than one
  processor),
\item Figure~\ref{peta_spdQRF} -- the best speedup obtained by
  \lstinline!PDGEQRF!, computed similarly, and
\item Figure~\ref{peta_cmp} -- the ratio of \lstinline!PDGEQRF!
  runtime to CAQR runtime.
\end{itemize}

As can be seen in Figure~\ref{peta_spdCAQR}, CAQR is expected to show
good scalability for large matrices.  For example, for $n = 10^{5.5}$,
a speedup of $1431$, measured with respect to the time on $2$
processors, is obtained on $8192$ processors. For $n=10^{6}$ a
speedup of $167$, measured with respect to the time on $32$
processors, is obtained on $8192$ processors.

In the technical report~\cite{TSQR_technical_report}, we also estimate
the fractions of time in computation, latency, and bandwidth for
\lstinline!PDGEQRF!  and CAQR.  These estimations show that for the
largest problems that can fit in memory, in the top left part of the
plots in Figure~\ref{peta_all}, the computation dominates the total
time, while in the right bottom part the latency dominates the total
time.  For the test cases situated between these two parts, the
bandwidth dominates the time.

CAQR leads to more significant improvements when the latency
represents an important fraction of the total time, the right bottom
part of Figure~\ref{peta_cmp}.  The best improvement is a factor of
$22.9$, obtained for $n = 10^4$ and $P = 8192$.  The speedup of the
best CAQR compared to the best \lstinline!PDGEQRF! for $n=10^4$ when
using at most $P=8192$ processors is larger than $8$, which is still
an important improvement.  The best performance of CAQR is obtained
for $P=4096$ processors and the best performance of
\lstinline!PDGEQRF! is obtained for $P=16$ processors.

Useful improvements are also obtained for larger matrices.  For $n =
10^6$, CAQR outperforms \lstinline!PDGEQRF! by a factor of $1.4$.
When the computation dominates the parallel time,
Figure~\ref{peta_cmp} predicts that there is no benefit from using
CAQR.  However, CAQR is never slower.  For any fixed $n$, we can take
the number of processors $P$ for which \lstinline!PDGEQRF!  would
perform the best, and measure the speedup of CAQR over
\lstinline!PDGEQRF! using that number of processors.  We do this in
Table \ref{tbl:CAQR:par:Peta:best}, which predicts that CAQR always is
at least as fast as \lstinline!PDGEQRF!, and often significantly
faster (up to $7.4 \times$ faster in some cases).

\begin{table}
  \begin{center}
\begin{tabular}{r|c|c}
$\log_{10} n$ & Best $\log_2 P$ for \lstinline!PDGEQRF! & CAQR speedup
\\ \hline
3.0 & 1     & 1 \\
3.5 & 2--3  & 1.1--1.5 \\
4.0 & 4--5  & 1.7--2.5 \\
4.5 & 7--10 & 2.7--6.6 \\
5.0 & 11--13 & 4.1--7.4 \\
5.5 & 13     & 3.0       \\
6.0 & 13     & 1.4       \\
\end{tabular}
\caption{Estimated runtime of \lstinline!PDGEQRF! divided by estimated
  runtime of parallel CAQR on a square $n \times n$ matrix, on the Peta platform,
  for those values of $P$ (number of processors) for which
  \lstinline!PDGEQRF! performs the best for that problem size.}
\label{tbl:CAQR:par:Peta:best}
  \end{center}
\end{table}

%% file: pr-conclusions.tex
\section{Conclusions and Future Work}
\label{pr:S:concl}

We have presented sequential and parallel algorithms that minimize the
communication performed during the QR factorization of tall and skinny
matrices and general rectangular matrices.  In the accompanying
paper~\cite{THEORY} we have shown that the new algorithms are optimal
in the amount of communication they perform, thus they are superior in
theory over existing algorithms.  In this paper we have presented
implementations demonstrating in practice significant speedups over
LAPACK and ScaLAPACK.  In particular, we have studied the performance
of parallel TSQR on a binary tree and sequential TSQR on a flat tree.

Implementations of parallel CAQR are currently underway.  Optimization
of the TSQR reduction tree for more general, practical architectures
(such as multicore, multisocket, or GPUs) is future work, as well as
optimization of the rest of CAQR to the most general architectures,
with proofs of optimality.

It is natural to ask to how much of dense linear algebra one can
extend the results of this paper, that is finding algorithms that
attain communication lower bounds.  In the case of parallel LU with
pivoting, refer to the paper by Grigori, Demmel, and
Xiang~\cite{grigori2008calu}, and in the case of sequential LU, refer
to the paper by Toledo \cite{toledo1997locality}.  More broadly, we
hope to extend the results of this paper to other linear algebra
operations, including two-sided factorizations (such as reduction to
symmetric tridiagonal, bidiagonal, or (generalized) upper Hessenberg
forms).  Once a matrix is symmetric tridiagonal (or bidiagonal) and so
takes little memory, fast algorithms for the eigenproblem (or SVD) are
available. Most challenging is likely to be finding eigenvalues of a
matrix in upper Hessenberg form (or of a matrix pencil).